# On noetherianity for logical formulas over fields

A. Berzins

University of Latvia

e-mail aberzins@latnet.lv

June 19, 2007

AMS Mathematics Subject Classification: 14A99, 14P05

Keywords: noetherianity, field

### Abstract

In this paper we consider noetherianity for formulas of propositional and predicate calculus over different fields. Three types of noetherianity are considered: standard noetherianity, logical noetherianity and denumerable noetherianity.

## 1 Introduction

Main concepts used in this article we define for arbitrary variety of algebras Θ. We shall consider different variants of noetherianity. Questions about different variants of noetherianity set B. Plotkin in his works [2, 3]. Because these concepts are new at first we must investigate classical case – variety of commutative algebras over field P. Main difference from classical noetherianity is that we consider not only equations but also logical formulas over an algebra from variety Θ. We see that question about some variants of noetherianity differs even for fields.

## 2 Definitions

We shall consider three types of formulas.

**1.** Equation in $\Theta$ is an equality $w = w'$ where $w, w' \in W(X)$; $W(X)$ is the free algebra in variety $\Theta$. Set of equations denote by $Eq(X)$; $X$ – finite set of variables. If $H$ is arbitrary algebra from $\Theta$, then we have affine space over $H$; it is $H^X$ or $Hom(W(X), H)$.

Every set of equations $T \subset Eq(X)$ determines affine set $T'_H = A$ -- the set of points in $H^X$, for which all equations from $T$ holds. Conversely, for every set $A \subset H^X$ denote by $A'$ the set of all equations which holds for every point from $A$. We get Galois correspondence between affine sets in $H^X$ and closed sets in $Eq(X)$.

**2.** In the second case in place of equations we consider formulas of propositional calculus over algebra $W(X)$. The set of all such formulas denote by $F_0(X)$. Every set of formulas $T \subset F_0(X)$ determines the subset of $H^X$ -- set of points for which all formulas from $T$ holds. Denote this set by $T_H^E$. The set $A$ which can be represented in form $A = T_H^E$ we shall call elementary set. Also every subset $A \subset H^X$ determines subset of formulas $A_H^E \subset F_0(X)$ which holds for all points of $A$. So we get Galois correspondence between elementary sets in $H^X$ and closed sets in $F_0(X)$.

**3.** In the third case we consider predicate calculus formulas. We consider all formulas over the free algebra $W(Y, X)$, where $Y$ is denumerable set of connected variables and $X$ is finite set of free variables. The set of all such formulas we denote by $F(Y, X)$. As in previous case every set of formulas $T \subset F(Y, X)$ define subset of $H^X$ -- the set of points for which all formulas from $T$ holds. Denote this set by $T_H^L$. The set $A$, which can be represented as $T_H^L$ we shall call simple set. Every subset $A \subset H^X$ determines a subset of formulas $A_H^L \subset F(Y, X)$ which holds for all points of $A$. So we get Galois correspondence between simple sets in $H^X$ and closed sets in $F(Y, X)$.

Also we have three variants of noetherianity.

**1.** NET – classical noetherianity for algebra $H$.

We say that algebra $H$ is noetherian if for every $T \subset Eq(X)$ $(T \subset F_0(X), T \subset F(Y,X))$ exist finite subset $T_0 \subset T$, such that $T_{0H}'' = T_H''$ ($T_{0H}^{EE} = T_H^{EE}$, $T_{0H}^{LL} = T_H^{LL}$). It means that for every congruence exist finite basis.

**2.** LOGNET – logical noetherianity for algebra $H$.

We say that that algebra $H$ is logical noetherianity if for every formula

$$\bigwedge_{i \in I} f_i \Rightarrow f \quad (f_i, f \in Eq(X), f_i, f \in F_0(X), f_i, f \in F(Y,X))$$

exist a finite subset $J \subset I$, such that $\bigwedge_{i \in J} f_i \Rightarrow f$.

$f \Rightarrow g$ means, that formula holds in $H^X$.

**3.** DENNET – denumerable noetherianity for algebra $H$.

It is the same as LOGNET, but the set $I$ is denumerable.

We wheel write $NET_0$, $LOGNET_0$, $DENNET_0$ for propositional calculus formulas.

## 3. Results

Example. Formula $\bigwedge_{a \in P} \neg (x = a) \Rightarrow (1 = 0)$ holds for every field $P$. However for every infinite field infinite conjunction we can not change on finite conjunction.

Corollary.

**Theorem 1.** For every infinite field $P$ does not hold LOGNET and $LOGNET_0$.

Now consider DENNET and $DENNET_0$.

**Theorem 2.** For field of real numbers $R$ does not hold DENNET.

Proof. We can not change implication $\bigwedge_{i \in N} \exists y (i(x^2 + y^2) = 1) \Rightarrow (x = 0)$ on finite conjunction.

**Definition.** We call subset $A \subset K^X$ pseudoaffine variety if it defines by formula. $\bigwedge_{i \in I} \varepsilon_i (f_i = 0)$, I finite, $\varepsilon_i$ -- empty or $\neg$.

Now consider question about denumerable noetherianity for nondenumerable fields. We take formula $\bigwedge_{i \in I} \varepsilon_i f_i \Rightarrow g$ and $f_i$ write in conjunctive form but $g$ in disjunctive form. We have

$$\bigwedge_{i \in N} (\bigwedge_k (\bigvee_l \varepsilon_{ikl}(f_{ikl} = 0))) \Rightarrow \bigvee_m (\bigwedge_n \varepsilon_{mn}(g_{mn} = 0))$$

Using formula $(a \Rightarrow b) \Leftrightarrow (\neg a \vee b)$ we change previous formula in form

$$\bigvee_{i \in N} (\bigvee_k (\bigwedge_l \neg \varepsilon_{ikl}(f_{ikl} = 0))) \vee (\bigvee_m (\bigwedge_n \varepsilon_{mn}(g_{mn} = 0)))$$

**Theorem 3.** Let $P$ be a nondenumerable field. If denumerable union of pseudoaffine varieties cover all affine space $P^X$, then we can take a finite subset of pseudoaffine varieties who cover all affine space.

Proof. Let $\bigcup_{i \in N} V_i = P^n$, $V_i$ -- pseudoaffine variety. Because number of $V_i$ is denumerable, then there exist $V_j$ which is everywhere dens in $P^n$ in Zarissky topology. We denote $V_j = V_1$. So $P^n \setminus V_1$ belongs to finite union of irreducible varieties of dimension $n-1$: $U_1, U_2, \ldots, U_k$. Every variety $U_i$ has some $V_{i_k}$, such that $U_i \cap V_{i_k}$ is everywhere dens in $U_i$. It means that $U_i \setminus V_{i_k}$ belongs to finite union of irreducible varieties of dimension $n-2$. We continue this process while take uncovered finite number of points, which we can cover by finite number of varieties $V_i$.

**Corollary.** For formulas of propositional calculus and nondenumerable fields denumerable noetherianity holds.

Follows from theorem 3 and transformation implication in denumerable disjunction.

And now we shall consider denumerable noetherianity for predicate calculus formulas. I note that for field $R$ denumerable noetherianity does not hold (theorem 2.). Next theorem show that for algebraic closed fields it holds.

**Theorem 4.** For algebraic closed fields and predicate calculus formulas denumerable noetherianity holds.

Proof. Follows from corollary to theorem 3 and next lemma (see [1]).

**Lemma.** For algebraic closed field every predicate calculus formula is equivalent to some propositional calculus formula.